\documentclass{amsart}

\newcommand{\ch}{\operatorname{ch}}
\newcommand{\coh}{cohomology}
\newcommand{\coin}{\operatorname{Coind}}
\newcommand{\coind}{\operatorname{Coind}^\gtg_\gth}
\newcommand{\Coker}{\operatorname{Coker}}
\newcommand{\Deg}{\operatorname{Deg}}
\newcommand{\gl}{\mathfrak{g}\mathfrak{l}_{\operatorname{res}}}
\newcommand{\glt}{\widetilde{\mathfrak{g}\mathfrak{l}}_{\operatorname{res}}}
\newcommand{\Hom}{\operatorname{Hom}}
\newcommand{\im}{\operatorname{Im}}
\newcommand{\ind}{\operatorname{Ind}^\gtg_\gth}
\newcommand{\inn}{\operatorname{Ind}}
\newcommand{\Ker}{\operatorname{Ker}}
\newcommand{\lb}{\mathcal{L}_\beta}
\newcommand{\mult}{\operatorname{mult}}
\newcommand{\op}{\operatorname{op}}
\newcommand{\rank}{\operatorname{rank}}
\newcommand{\sh}{H^{\infty/2+ \bullet}}
\newcommand{\si}{semi-in\-fi\-nite}
\newcommand{\sindh}{\operatorname{S-ind}^{\gtg}_{\gth}}
\newcommand{\sind}{\operatorname{S-ind}}
\newcommand{\sj}{semi\-jec\-tive}

\newcommand{\US}{\operatorname{US}}
\newcommand{\uea}{universal enveloping algebra}

\newcommand{\gta}{\mathfrak{a}}
\newcommand{\gtg}{\mathfrak{g}}
\newcommand{\gth}{\mathfrak{h}}
\newcommand{\gtn}{\mathfrak{n}}
\newcommand{\gtt}{\mathfrak{t}}
\newcommand{\nc}{\mathbb{C}}
\newcommand{\nz}{\mathbb{Z}}
\newcommand{\OO}{\mathcal{O}_0}
\newcommand{\wgtg}{\widehat{\mathfrak{g}}}

\font\cyr=wncyi8

\newtheorem{thm}{Theorem}[section]
\newtheorem{prop}[thm]{Proposition}

\theoremstyle{definition}
\newtheorem{df}{Definition}[section]
\newtheorem{ex}{Example}[section]

\theoremstyle{remark}
\newtheorem*{rem}{Remark}
\newtheorem*{ack}{Acknowledgments}

\numberwithin{equation}{section}

\begin{document}

\title{Semi-Infinite Induction and Wakimoto Modules}
\author{Alexander A. Voronov}
\thanks{Research supported in part by an AMS Centennial Fellowship and 
NSF grant DMS-9402076.}
\address{Department of Mathematics \\
Michigan State University\\
East Lansing, MI 48824-1027}
\email{voronov@math.msu.edu}
\urladdr{http://www.mth.msu.edu/$\,\tilde{}\,$voronov/}

\begin{abstract}
The purpose of this paper is to suggest the construction and study
properties of semi-infinite induction, which relates to semi-infinite
cohomology the same way induction relates to homology and coinduction
to cohomology. We prove a version of the Shapiro Lemma, showing that the
semi-infinite cohomology of a module is isomorphic to that of the
semi-infinitely
induced module. A practical outcome of our construction is a simple
construction of the Wakimoto modules, highest-weight modules used in
double-sided BGG resolutions of irreducible modules.
\end{abstract}

\maketitle

Semi-infinite \coh\ of Lie algebras, introduced as the appropriate
mathematical setting for BRST theory by B.~L. Feigin \cite{feigin}
(see also \cite{fgz,ya} regarding basic facts on \si\ \coh), is a
\coh\ theory that has properties in common with both \coh\ and homology. 
Semi-infinite \coh\ has become an important tool in representation
theory of Lie algebras and quantum groups and string theory, see, for
instance,
\cite{ar:bgg,ar:quantum,ar:quantumII,bv,ffr,jose2,lz:appl,soergel}.

The purpose of this paper is to suggest the construction and study
properties of \si\ induction, which relates to \si\ \coh\ the same way
induction relates to homology and coinduction to \coh. The proof of
our main theorem (Theorem~\ref{main}, the \si\ Shapiro Lemma) is based
on the independence of the choice of a resolution, which follows from
the machinery of \si\ homological algebra developed in \cite{ya}. A
practical outcome of our construction is a simple construction of
Wakimoto modules, which were constructed by Feigin and E.~Frenkel
\cite{ffr} in rather roundabout terms: using bosonization and also as
$H^{\infty/2+0}_{U} (X, \mathcal L_\lambda)$, a hypothetical \si\
\coh, with support on the big Schubert cell, of an invertible sheaf
over a \si\ flag manifold. The idea that Wakimoto modules might be
obtained by some kind of \si\ induction goes back to the original
paper of Feigin and Frenkel \cite{ffr}. This idea was implemented by
S.~M. Arkhipov \cite{ar:bgg}, who suggested an indirect construction
of Wakimoto modules, which, in fact, may be considered as representing
a different approach to \si\ induction. Wakimoto modules play an
intermediate role between Verma and contragredient Verma modules: they
all have a very similar behavior with respect to \si\ \coh, usual
homology, and \coh, respectively. The three types of modules have the
same character (see Proposition~\ref{structure}) but a different
layout of irreducible pieces.

This work was motivated in part by a construction of N.~Berkovits and
C.~Vafa \cite{bv} of $N=1$ and $N=2$ string theories out of a given
$N=0$ (bosonic) string theory. In that construction, the bosonic
string arose as a particular class of vacua for the $N=1$ string and
the $N=1$ as a particular class of vacua for the $N=2$ string.
Berkovits and Vafa also suggested that there must be a universal
string theory comprising all possible string theories, including, for
instance, $W_N$ strings, as particular choices of vacua. It was
J.~Figueroa-O'Farrill \cite{jose2,jose} who observed that
Berkovits-Vafa's construction presumably had to do with some sort of
\si\ induction and inquired whether such construction had been
known. This paper contains an answer to his question. However, the
physically oriented reader must be warned that the physics problem is
more complex than the mathematical model. The physical spaces of vacua
of two string theories in the hierarchy of the conjectural universal
string theory should not only match as vector spaces, but the
correlators between different vacua should be equal to each other ---
this would guarantee that the two theories are physically the same.
Mathematically, this amounts to the problem of performing \si\
induction in the category of vertex operator algebras: assume that the
module over a smaller algebra is a vertex operator algebra and
construct a matching vertex operator algebra structure on the
semi-induced module. We do not know how to do that.

Throughout the paper we will be working over the ground field $\nc$ of
complex numbers.

\begin{ack}
We are very grateful to Jos\'e Figueroa-O'Farrill, who brought the
question to our attention, Victor Ginzburg of Chicago, who pointed out
the papers of S.~M. Arkhipov and W.~Soergel, S.~M. Arkhipov and
T.~Kimura, who made very helpful remarks, A.~S. Schwarz, who
apparently convinced us to write down the results of this paper,
W.~Soergel, who spotted some inaccuracies in the proofs of two
propositions of \cite{ya}, which have been corrected in the present
paper, see Propositions \ref{iso} and \ref{iso-1}, and A.~Sevastianov,
who pointed out an error in the original version of the proof of
Semi-infinite Shapiro Lemma.  We would also like to thank the Institut
des Hautes \'Etudes Scientifiques at Bures-sur-Yvette for its
hospitality during the summer of 1996, when the essential part of the
work on this paper was done.
\end{ack}

\section{Semi-infinite induction}

Induction is a certain construction of ``base change'' in
representation theory. Given a Lie algebra and a subalgebra of it,
$\gth \subset \gtg$, as well as an $\gth$-module $M$, one constructs a
$\gtg$-module $\ind M$, such that the corresponding homology does not
change: $H_\bullet (\gtg, \ind M) = H_\bullet (\gth, M)$.  Coinduction
starts from the same data, but constructs a $\gtg$-module $\coind M$
possessing the similar property with respect to cohomology: $H^\bullet
(\gtg, \coind M) = H^\bullet (\gth, M)$.  Semi-infinite induction will
do the same, pertinent to \si\ \coh: $\sh (\gtg, \sindh M) = \sh
(\gth, M)$.

\subsection{Coinduction}
Let us briefly recall the coinduction construction. The {\it coinduced
module} is defined as
\[
\coind M = \Hom_\gth (U(\gtg), M),
\]
where $\gtg$ acts on the \uea\ $U(\gtg)$ by multiplication on the
right and $\gth$ on the left. The following theorem is a standard fact
of homological algebra.

\begin{thm}[Shapiro Lemma]
There is a natural isomorphism of cohomology
\[
H^\bullet (\gtg, \coind M) = H^\bullet (\gth, M).
\]
\end{thm}

\begin{proof}
Take a projective resolution
\[
P^\bullet: \qquad \dots \to P^{-1} \to P^0 \to 0
\]
of the trivial $\gtg$-module $\nc$. It is a bounded above complex of
projective $\gtg$-modules such that $H^\bullet (P^\bullet) = H^0
(P^\bullet) = \nc$. Notice that it is also a projective resolution of
$\nc$ in the category of $\gth$-modules and $\Hom (P^\bullet, M)$ is
an injective resolution of the $\gth$-module $M$. Therefore, one
produces the cohomology of $\gth$ with coefficients in $M$ by applying
the functor of $\gth$-invariants to this resolution and then taking
the \coh\ of the obtained complex:
\[
H^\bullet (\gth, M) = H^\bullet (\Hom_{\gth} (P^\bullet,M)).
\]

On the other hand, $\Hom (P^\bullet, \coind M)$ is an injective
resolution of the $\gtg$-module $\coind M$. Thus, if we apply the
functor of $\gtg$-invariants, we will get the \coh\ of $\gtg$ with
coefficients in $\coind M$:
\[
H^\bullet (\gtg, \coind M) = H^\bullet (\Hom_{\gtg} (P^\bullet, \coind
M)).
\]

It remains to notice that the two complexes are naturally isomorphic
because of the universality property of coinduction:
\begin{equation}
\label{universal}
\Hom_{\gtg} (P^\bullet, \coind M) = \Hom_{\gth} (P^\bullet,M).
\end{equation}
\end{proof}

\subsection{Semi-infinite structure and semi-invariants}
\label{si-str}

The \si\ induction and Shapiro Lemma are very much parallel to the
classical coinduction case, except that each step calls for an
entirely new ingredient. The construction of a semi-induced module
requires an intermediate object between the \uea\ and its dual. The
\si\ analogue of Shapiro Lemma needs a functor, like the one of
invariants, to get back to the original module from an induced one.
The proof of it also needs the machinery of two-sided resolutions and
two-sided derived functors. All these ingredients have been developed
in the earlier paper \cite{ya}.

Suppose we have a Lie algebra with a {\it \si\ structure}. This may be
understood as a graded Lie algebra $\gtg = \bigoplus_{n \in \nz}
\gtg_n$ with finite-dimensional graded components, along with the
following structure. Decompose the algebra $\gtg$ into the direct sum
of two subalgebras:
\begin{gather*}
\gtg = \gtg_+ \oplus \gtg_-,\\
\gtg_+ = \bigoplus_{n > 0} \gtg_n, \quad
\gtg_- = \bigoplus_{n \le 0} \gtg_n.
\end{gather*}
Assume that the natural mapping $\gtg \to \gl$ via the adjoint
representation is lifted to a mapping $\gtg \to \glt$, where $\gl =
\gl (\gtg)$ is a ``restricted'' general linear algebra of the vector
space $\gtg$, for example, the one consisting of
operators whose $\gtg_+ \to \gtg_-$ block is of a finite rank; $\glt$
is a nontrivial central extension of $\gl$, see e.g., \cite{adkp}. If
this lifting is not possible, we should replace $\gtg$ by the
corresponding central extension, which will be lifted canonically. Let
$\beta: \gtg \to \nc$ be the linear functional defined by this lifting
and a splitting of the extension $0 \to \nc \to \glt \to \gl \to 0$ as
an extension of vector spaces. In fact, this splitting can be chosen
in such a way that $\beta$ vanishes on all $\gtg_n$ but $\gtg_0$, see
for example, \cite[Proposition 2.4]{ya}. We will assume this for the
sake of simplicity.  Notice that $\beta$ defines a one-cocycle on
$\gtg_-$ and the zero one-cocycle on $\gtg_+$. Denote the
corresponding one-dimensional modules by one symbol $\lb$. As vector
spaces, the modules $\lb$ are canonically isomorphic to $\nc$. We
assume the functional $\beta$ to be part of a \si\ structure on a Lie
algebra.

Throughout this paper, a \emph{$\nz$-graded vector space} means a
vector space $M$ with a collection of subspaces $M_n$, $n \in \nz$,
such that $\bigoplus_{n \in \nz} M_n \subset M \subset \prod_{n \in
\nz} M_n$. For two graded vector spaces $M$ and $N$, we will define
the space $\Hom(M,N)$ as follows:
\begin{equation}
\label{hom}
\Hom (M,N) = \bigoplus_{n > 0 } \Hom(M_n, N) \oplus \prod_{n\le 0}
\Hom(M_n,N).
\end{equation}
It has a natural $\nz$-bigrading. This definition of $\Hom$ is
motivated by the following argument. If we consider the inverse-limit
topology on $M$, that is, the topology coming from $\varprojlim_{n \to
\infty} M/F^n M$, where $F^n M = M \cap (\prod_{k\ge n} M_k)$, and the
discrete topology on $N$, the space of continuous linear maps $M \to
N$ will be identified with $\Hom(M,N)$ in \eqref{hom}.

All $\Hom$'s and duals will be understood in this sense throughout the
paper. Note that with this definition of $\Hom$, for $\nz$-graded
vector spaces $A$, $S$, and $M$, such that $A = \bigoplus_{n \le N}
A_n$ and $M = \bigoplus_{n \le N'} M_n$ and have finite-dimensional
graded components $A_n$ and $M_n$ for all $n$, one has canonical
isomorphisms
\begin{gather}
\Hom (S \otimes M, A) \xrightarrow{\sim} \Hom (S, \Hom (M,A)),
\nonumber \\
\label{conj}
\Hom (A, S \otimes M) \xrightarrow{\sim} \Hom (M^* \otimes A, S)
\xleftarrow{\sim} \Hom(\Hom (M,A), S).
\end{gather}

We will consider $\nz$-graded $\gtg$-modules, as well as objects of
the so-called \emph{category $\OO$}, by which we will here mean the
category of $\nz$-graded $\gtg$-modules $M$ semisimple over $\gtg_0$,
such that they are direct sums of finite-dimensional graded components
$M_n$ and the grading is bounded above: $n \le N$:
\[
M = \bigoplus_{n=-\infty}^N M_n .
\]


Given a \si\ structure on a Lie algebra $\gtg$, we can define the
functor of \emph{semi-invariants} of a $\gtg$-module $M$:
\begin{equation}
\label{sinv}
\left(M \right)^{\gtg_+}_{\gtg_-} = \im \left( (M \otimes \lb)^{\gtg_+}
\to (M \otimes \lb)_{\gtg_-} \right),
\end{equation}
the image of the natural projection of the $\gtg_+$-invariants $\{x
\in M \otimes \lb \; | \; gx = 0 \text{ for all } g \in \gtg_+\}$ onto the
$\gtg_-$-coinvariants $(M \otimes \lb)_{\gtg_-} = M \otimes
\lb/{\gtg_-} (M \otimes \lb)$. This functor coincides with 
the \si\ \coh\ group $H^{\infty/2 +0} (\gtg, M)$, see
Section~\ref{si-coh}, for a class of modules $M$ called
\emph{\sj} in \cite[Corollary 3.4]{ya}. We will recall this 
notion along with a more general concept of a \sj\ resolution below in
Section~\ref{semijective}.

\subsection{Semi-infinite \coh}
\label{si-coh}

Suppose $\gtg$ is a Lie algebra with a \si\ structure. The space
$\Lambda^{\infty/2+\bullet} \gtg$ of \emph{\si\ forms} on $\gtg$ is an
irreducible representation of the Clifford algebra based on the vector
space $\gtg \oplus \gtg^*$, where $\gtg^* = \Hom (\gtg, \nc)$ in the
sense of \eqref{hom}. This representation is spanned by a vector
$\omega_0$, called the vacuum, satisfying the conditions
\begin{align*}
g^* \omega_0 & = 0 && \text{ for $g^* \in \gtg_+^*$,}\\ g \omega_0 & =
0 && \text{ for $g \in \gtg_-$.}
\end{align*}
The space of \si\ forms is graded: $\Lambda^{\infty/2+\bullet} \gtg =
\bigoplus_p \Lambda^{\infty/2+p} \gtg$, the degree being counted from
$\Deg \omega_0 = 0$, $\Deg g = 1$ for $g \in \gtg$ and $\Deg g^* = -1$
for $g^* \in \gtg^*$. One can choose a basis $\{e_i \; | \; i \in
\nz\}$ in $\gtg$ compatible with the $\nz$-grading and think of the
vacuum as $\omega_0 = e_0 \wedge e_{-1} \wedge e_{-2} \wedge \dots$
Then an arbitrary \si\ form will be a sum of forms $\omega = g_1
\wedge g_2 \wedge \dots$, where $g_i \in \gtg$, and $\omega$ and
$\omega_0$ have equal \si\ tails, that is to say, terms in the wedge
products coincide starting from some point on. The choice of a basis
determines an inner product on $\gtg$ with respect to which the basis
is orthonormal.

Define the \emph{standard \si\ complex} for our Lie algebra $\gtg$
with coefficients in a $\nz$-graded $\gtg$-module $M$. This is the
complex
\[
C^{\infty/2+n} = \Hom(\Lambda^{\infty/2+n} \gtg, M), \qquad n \in \nz,
\]
where the ``$\Hom$'' is understood as in \eqref{hom} with respect to
the interior grading on $\Lambda^{\infty/2+\bullet} \gtg$ given by
$\deg \omega_0 = 0$, $\deg g = i$ for $g \in \gtg_i$, and $\deg g^* =
-i$ for $g^* \in \gtg_i^* := (\gtg_i)^*$. The differential is defined
as
\begin{multline*}
(d \phi) (g_1 \wedge g_2 \wedge \dots) \\ = \sum_i (-1)^{i+1} (g_i
+\beta(g_i)) \phi(g_1\wedge g_2 \wedge \dots \wedge \widehat g_i
\wedge
\dots)
\\+ \sum_{i < j} (-1)^{i+j} \phi(\, : \! [g_i,g_j] \wedge g_1 \wedge
\dots
\wedge \widehat g_i \wedge \dots \wedge \widehat g_j \wedge \dots:),
\end{multline*}
where $g_i$'s are assumed to be homogeneous, and the normal ordering
sign $::$ means that whenever $[g_i,g_j] \in \gtg_-$
\[
:[g_i,g_j] \wedge g_1 \wedge \dots \wedge \widehat g_i \wedge \dots
\wedge \widehat g_j \wedge \dots: = P_{ij}([g_i,g_j]) \wedge g_1 \wedge
\dots \wedge \widehat g_i \wedge \dots
\wedge \widehat g_j \wedge \dots,
\]
$P_{ij}$ being the orthogonal projection onto the orthogonal
complement to the linear span of $g_i$ and $g_j$ in the sense of the
bilinear inner product determined by the choice of a basis in $\gtg$.
Obviously, either sum in $d\phi$ will be finite for each $\phi$. Note
that $M$ needs not lie in $\OO$ or even have a bounded above
$\nz$-grading in this definition. This relaxes the hypothesis usually
made in the literature, see \cite{feigin,fgz}.

\subsection{Universal \sj\ module}
\label{US}

We are going to define a bimodule which plays the same role with
respect to \si\ \coh\ as the \uea\ plays with respect to homology.
Suppose $\gtg$ is a Lie algebra with a \si\ structure. The following
$U(\gtg)$-bimodule will be called \emph{universal semijective}:
\[
\US = \left(\Hom (U(\gtg), U(\gtg))\right)_{\gtg_+}^{\gtg_-},
\]
where the ``anti'' $\gtg$-semi-invariants are taken with respect to
the action
\begin{equation}
\label{g-act}
(g f) (u) = g (f(u)) + f (ug) \qquad \text{for $f \in \Hom (U(\gtg),
U(\gtg))$.}
\end{equation}
Notice that the ``anti'' semi-invariants with respect to $\gtg$ are
understood as in \eqref{sinv}, the functional $\beta$ being replaced by
$-\beta$.

The $U(\gtg)$-bimodule structure on $\US$ is defined through the
natural left $\gtg$-action
\begin{equation}
\label{left}
(g f) (u) = f(-g u) \qquad \text{for $f \in \Hom (U(\gtg), U(\gtg))$}
\end{equation}
and the natural right $\gtg$-action
\begin{equation}
\label{right}
(fg) (u) = f(u)g \qquad \text{for $f \in \Hom (U(\gtg), U(\gtg))$.}
\end{equation}

The universal semijective module $\US$ was introduced in \cite{ya}
under the name ``standard \sj'' and denoted $\operatorname{SS}$.
Arkhipov \cite{ar:bgg} and Soergel \cite{soergel}, who also suggested
different constructions for it, used it under the names of
\emph{semiregular module} and \emph{semi-regul\"are Bimodul},
respectively. The following two important statements regarding the
structure of the left and the right actions on $\US$ were given
inaccurate proofs in \cite{ya}, as was pointed out by Soergel; see
also a different proof of the composite isomorphism
$U(\gtg_+)^*\otimes_{\gtg_+} U(\gtg) =
\Hom_{\gtg_-}(U(\gtg),U(\gtg_-)\otimes \mathcal{L}_{-\beta})$ of
Propositions \ref{iso} and \ref{iso-1} under stronger assumptions on
$\gtg$ in Soergel \cite{soergel}.

\begin{prop}
\label{iso}
As a right $\gtg$-module, the universal semijective module $\US$ is
isomorphic to $U(\gtg_+)^*\otimes_{\gtg_+} U(\gtg)$, where $\gtg_+$
acts as in \eqref{g-act} and $\gtg$ acts as in \eqref{right}.
\end{prop}

\begin{proof}
Let us shorten the notation $U(\gtg)$ to $U$ and consider the standard
\si\ complex $C^{\infty/2+\bullet}_{\op}$ for the Lie algebra $\gtg$ with
the opposite \si\ structure (the opposite $\nz$-grading and $-\beta$
as the structure one-cochain) with coefficients in the $\gtg$-module
$M = \Hom(U,U)$, see \eqref{g-act}. The module $M$ is clearly \sj\
with respect to this action, and therefore (Corollary 3.4 of
\cite{ya}), its \si\ \coh\ is computed as the semi-invariants
$M^{\gtg_-}_{\gtg_+}= \US$ in degree 0 and 0 in all other degrees. We
will compute the same cohomology group as $U(\gtg_+)^*
\otimes_{\gtg_+} U$, using a spectral sequence. This will prove the
proposition.

This spectral sequence, a \si\ version of the Hochschild-Serre
spectral sequence for the pair $\gtg_- \subset \gtg$, was considered
in \cite[Theorem 2.2]{ya} and comes from the following filtration on
the standard \si\ complex $C^{\infty/2+\bullet}_{\op} \linebreak[1] =
\linebreak[0] \Hom(\Lambda^{\infty /2 + \bullet}_{\op} \gtg, M)$,
where $\op$ in the subscript reminds us to consider the opposite
grading on $\gtg$:
\begin{multline*}
F^p C^{\infty/2 + p + q}_{\op} \\:= \{ \phi \in C^{\infty/2 + p +
q}_{\op}
\; | \; b_1 \dots b_{q+1} \phi = 0 \text{ for each $b_1, \dots, b_{q+1}
\in \gtg_-$}\},
\end{multline*}
where $p \le 1$, $q \ge 0$ and the multiplication $b \phi$ is a
contraction of a \si\ form.  It is evident that
\[
C^{\infty/2 +\bullet}_{\op} \supset \dots \supset F^p \supset F^{p+1}
\supset \dots
\supset F^1 = 0,
\]
that is, the filtration is decreasing and regular, and moreover,
\[
\bigcup_p F^p = C^{\infty/2 +\bullet}_{\op}.
\]
One can check that the differential $d$ in $C^{\infty/2
+\bullet}_{\op}$ is compatible with the filtration. We have a filtered
complex, which gives rise to a spectral sequence $ E^{p,q}_r $
converging to $\sh_{\op} (\gtg;M) = M^{\gtg_-}_{\gtg_+}$. The first
term of this spectral sequence is computed as $E^{p,q}_1 = H^q
(\gtg_-, \,
\Hom(\Lambda_{\op}^{\infty/2+p} (\gtg/\gtg_-), M))$. Since the module $M$ is
$\gtg_-$-injective, $E^{p,q}_1 = 0$ for any $q > 0$. For $q=0$, we
have the row $E^{\bullet,0}_1$:
\[
\dots \to \Hom_{\gtg_-} (\Lambda_{\op}^{\infty/2-1} (\gtg/\gtg_-), M) \to
\Hom_{\gtg_-} (\Lambda^{\infty/2+0}_{\op} (\gtg/\gtg_-), M) \to 0.
\]
This complex can be identified with
\[
\dots \to (\gtg/\gtg_- \otimes \mathcal{L}_{-\beta} \otimes M)^{\gtg_-} \to
(\mathcal L_{-\beta} \otimes M)^{\gtg_-} \to 0,
\]
the differential being the $\gtg_+$-homology differential. Since
$M=\Hom(U,U)$, this complex can further be identified with
\[
\dots \to  \gtg_+ \otimes \Hom(U(\gtg_+), U) \to
\Hom(U(\gtg_+),U) \to 0,
\]
with the differential being the standard differential in the homology
complex of $\gtg_+$ with coefficients in $\Hom(U(\gtg_+),U)$.  This
shows that the 0th homology $E^{0,0}_2$ of this complex is
$\Hom(U(\gtg_+),U)_{\gtg_+} = U(\gtg_+)^* \otimes_{\gtg_+} U$. Since
$U(\gtg_+)^* \otimes U$ is free over $\gtg_+$, the higher homology
vanishes, i.e., $E^{p,0}_2=0$ for $p>0$ and, therefore, the
spectral sequence collapses at $E_2 = E_\infty$. Thus $U(\gtg_+)^*
\otimes_{\gtg_+} U = E_2^{0,0} = H^{\infty/2+0}(\gtg, M) =
M^{\gtg_-}_{\gtg_+}$. Note that all the identifications made do not
affect the right action \eqref{right} of $\gtg$, and, thus, we have an
isomorphism of $\gtg$-modules.
\end{proof}

\begin{prop}
\label{iso-1}
As a left $\gtg$-module, the universal semijective module $\US$ is
isomorphic to $\Hom_{\gtg_-}(U(\gtg),U(\gtg_-)\otimes \mathcal
L_{-\beta})$, where $\gtg_-$ acts as in \eqref{g-act} and $\gtg$ acts
as in
\eqref{left}.
\end{prop}

\begin{proof}
\begin{sloppypar}
The proof of this proposition is very similar to that of
Proposition~\ref{iso}. For the module $M = \Hom(U,U)$, we will compute
the \si\ \coh\ $M^{\gtg_-}_{\gtg_+}= \US$ differently, using a \si\
version of the Hochschild-Serre spectral sequence for $\gtg_+ \subset
\gtg$, see
\cite[Theorem 2.3]{ya}. This spectral sequence comes from the
following filtration on the standard \si\ complex
$C^{\infty/2+\bullet}_{\op}$:
\begin{multline*}
{F'}^p C^{\infty/2 + p + q}_{\op} \\:= \{ b_1^* \dots b_p^* \phi \in
C^{\infty/2 + p + q}_{\op} \; | \; b_1^*, \dots, b^*_{p} \in \gtg_-^*,
\; \phi \in C^{\infty/2 + q}_{\op}\},
\end{multline*}
where $p \ge 0$, $q \le 0$ and $b^* \phi$ is the exterior
multiplication. We have
\[
C^{\infty/2 +\bullet}_{\op} = {F'}^0 \supset \dots \supset {F'}^p
\supset {F'}^{p+1} \supset \dots \supset 0,
\]
which in particular means that the filtration is decreasing and
coregular; moreover,
\[
\bigcap_p {F'}^p = 0.
\]
This ensures convergence of the associated spectral sequence to
$E'_\infty = \linebreak[0] \sh_{\op}(\gtg; M) \linebreak[1] =
M^{\gtg_-}_{\gtg_+}$. The first term of this spectral sequence is
${E'}^{p,q}_1 = H^{\infty/2+q}_{\op} (\gtg_+, \, \Hom(\Lambda^{p}
(\gtg/\gtg_+), M)) = H_{-q} (\gtg_+, \,
\Hom(\Lambda^{p} (\gtg/\gtg_+), M) \otimes \mathcal L_{-\beta})$. Since the
module $M$ is $\gtg_+$-projective, ${E'}^{p,q}_1 = 0$ for any $q < 0$.
The $q=0$ row ${E'}^{\bullet,0}_1$ is then
\[
0 \to (M \otimes \mathcal L_{-\beta})_{\gtg_+} \to (\Hom (\gtg/\gtg_+,
M)\otimes \mathcal L_{-\beta})_{\gtg_+} \to \dots,
\]
the differential being the $\gtg_-$-\coh\ differential. Recalling that
$M=\Hom(U,U)$, we identify this complex with
\[
0 \to \Hom(U, U(\gtg_-) \otimes \mathcal L_{-\beta}) \to \Hom(\gtg_-, \Hom(U,
U(\gtg_-) \otimes \mathcal L_{-\beta})) \to \dots
\]
with the standard differential in the \coh\ complex of $\gtg_-$ with
coefficients in $\Hom(U,U(\gtg_-) \otimes \mathcal L_{-\beta})$.  Thus
the 0th \coh\ ${E'}^{0,0}_2$ of this complex is $\Hom_{\gtg_-}
(U,U(\gtg_-) \otimes \mathcal L_{-\beta})$. The higher \coh\
${E'}^{p,0}_2$ for $p>0$ vanishes, because the $\Hom$ module is
cofree over $\gtg_-$, and, therefore, the spectral sequence collapses
at ${E'}_2 = {E'}_\infty$. Thus $\Hom_{\gtg_-} (U, U(\gtg_-) \otimes
\mathcal L_{-\beta}) = {E'}_2^{0,0} = H^{\infty/2+0}(\gtg, M) =
M^{\gtg_-}_{\gtg_+}$. As in the proof of the previous proposition, all
the identifications made do not affect the left action \eqref{left} of
$\gtg$, and, thus, we have an isomorphism of $\gtg$-modules.
\end{sloppypar}
\end{proof}

Thus the universal semijective module $\US$ is free over $\gtg_-$ and
cofree over $\gtg_+$ with respect to either action of $\gtg$. Its
natural $\nz$-grading is also bounded above, and when $\gtg_0=0$,
$\US$ lies in $\OO$.

\subsection{Semi-infinite induction}

Suppose now we have a graded subalgebra $\gth$ in $\gtg$. The
subalgebra inherits a natural \si\ structure from the one on $\gtg$
with $\gth_+ = \gtg_+ \cap \gth$ and $\gth_- = \gtg_- \cap \gth$, so
that $\gth = \gth_+ \oplus \gth_-$.  We also assume that the
restriction of the one-cochain $\beta$ to $\gth$ defines a morphism
$\gth \to \glt (\gth)$. Let $M$ be an $\gth$-module also provided with
a \si\ structure, see Section~\ref{si-str}. Then the
\emph{semi-induced module} is defined as
\begin{equation*}
\sindh M := \US \otimes^{\gth_+}_{\gth_-} M,
\end{equation*}
where $\otimes^{\gth_+}_{\gth_-}$ means taking the semi-invariants of
the action of $\gth$ by
\begin{equation}
\label{h-action}
h (u\otimes m) = -uh \otimes m + u \otimes hm \qquad \text{on $\US
\otimes M$.}
\end{equation}
The space $\sindh M$ inherits the structure of a $\gtg$-module from
the left $\gtg$-action \eqref{left} on $\US$. It is easy to see that
if $M \in \OO$ as an $\gth$-module and $\gth_0 = \gtg_0$, then the
semi-induced module $\sindh M$ will also lie in the category $\OO$
with respect to the bigger algebra $\gtg$. Otherwise, its
$\nz$-grading will be just bounded above.

Consider a few degenerations of \si\ induction.
\begin{ex}[Coinduction]
Suppose $\gtg = \gtg_+$ and therefore $\gth = \gth_+$. Then $\US =
U(\gtg)^*$ and $\sindh M = U(\gtg)^* \otimes^{\gth} M = \coind M$.
\end{ex}

\begin{ex}[Induction]
Now suppose $\gtg = \gtg_-$ and $\gth = \gth_-$. Then $\US = U(\gtg)$
and $\sindh M = U(\gtg) \otimes_{\gth} M = \ind M$.
\end{ex}

\begin{ex}[The universal semijective module]
Let $\gth = 0$ and $M = \nc$. Then $\sind^\gtg_0 \nc = \US$.
\end{ex}

\subsection{Semijective resolutions and modules}
\label{semijective}

Here we are going to recall the notion of a \sj\ resolution from
\cite{ya}. By a \emph{resolution} of a module $M$, we as usual mean a
complex $S^\bullet$ whose cohomology $H^\bullet (S^\bullet)$ is
identified with $M[0]$, which denotes the complex $\dots \to 0 \to M
\to 0 \to \dots$ with $M$ placed in degree 0. A
\emph{\sj\ complex} is a complex $S^\bullet$ of $\gtg$-modules, such
that
\begin{enumerate}
\item it is $K$-\emph{injective} as a complex of $\gtg_+$-modules, 
\emph{i.e}., for every acyclic complex $A^\bullet$ of $\gtg_+$-modules,
$\Hom_{K(\gtg_+)} (A^\bullet,\, S^\bullet) = 0$, where $K(\gtg_+)$ is
the homotopy category of complexes of $\gtg_+$-modules;
\item it is $K$-\emph{projective relative to} $\gtg_+$, \emph{i.e}.,
for every acyclic complex $B^\bullet$ of $\gtg$-modules which is
isomorphic to $0$ in the category $K(\gtg_+)$, $\Hom_{K(\gtg)}
(S^\bullet,\, B^\bullet) = 0$.
\end{enumerate}
A \emph{\sj\ module} is nothing but a \sj\ complex $0 \to M \to
0$. This is equivalent to being injective as a $\gtg_+$-module and
projective relative to $\gtg_+$ as a $\gtg$-module. The universal
semijective module $\US$ is an example, with respect to either action
of $\gtg$, see Section~\ref{US}. A serious difficulty in dealing with
unbounded complexes comes from the fact that a complex made up of \sj\
modules will not necessarily be \sj. This difficulty can be gotten
round by the above notion of a \sj\ complex.

Define a \emph{weakly \sj\ complex} as a complex of $\gtg$-modules
which is $K$-injective as a complex of $\gtg_+$-modules and
$K$-projective as a complex of $\gtg_-$-modules. It is straightforward
to see that every \sj\ complex is weakly \sj. Note also that the
functor of semi-invariants takes acyclic weakly semijective complexes
to acyclic ones. The proof of this statement is similar to the proof
of vanishing theorem \cite[Theorem 2.1]{ya} for \si\ \coh\ of a weakly
\sj\ module, except that now one applies the functor of \si\ \coh\ to
an acyclic weakly \sj\ complex. Therefore, \si\ \coh\ may in fact be
defined by applying the functor of semi-invariants to weakly \sj\
resolutions termwise. We will use these facts in the proof of
Semi-infinite Shapiro Lemma in the following section.

\subsection{Semi-infinite Shapiro Lemma}

In the following theorem, we will assume that $\gtg_0 = \gth_0$ to
make sure that $\sindh M \in \OO$ whenever $M \in \OO$.

\begin{thm}[Semi-infinite Shapiro Lemma]
\label{main}
For any $\gth$-module $M$ in the category $\OO$, there exists a
canonical isomorphism
\[
\sh (\gtg,\sindh M) = \sh (\gth, M).
\]
\end{thm}

\begin{proof}
Take a weakly \sj\ resolution $S^\bullet = \dots \to S^{-1} \to S^0
\to S^1 \to \dots$ of the trivial $\gtg$-module $\nc$. Notice that
$S^\bullet$ is also a weakly \sj\ resolution of $\nc$ in the category
of $\gth$-modules. Indeed, for an acyclic complex $A^\bullet$ of
$\gth_+$-modules, the complex $U(\gtg_+) \otimes_{\gth_+} A^\bullet$
of $\gtg_+$-modules will also be acyclic, and $\Hom_{K(\gth_+)}
(A^\bullet, S^\bullet) \linebreak[1] = \linebreak[0] \Hom_{K(\gtg_+)}
(U(\gtg_+) \otimes_{\gth_+} A^\bullet, S^\bullet) = 0$, because
$S^\bullet$ is $K$-injective as a complex of $\gtg_+$-modules. This
shows that $S^\bullet$ is $K$-injective as a complex of
$\gth_+$-modules. Similarly, one proves that $S^\bullet$ is
$K$-projective as a complex of $\gth_-$-modules.

Now notice that $S^\bullet \otimes M$ will be a weakly \sj\ resolution
of $M$ over $\gth$. To see that $S^\bullet \otimes M$ is a weakly \sj\
complex of $\gth$-modules, we have to check two things, as in the
previous paragraph. The first is that $\Hom_{K(\gth_+)} (A^\bullet ,
S^\bullet \otimes M) = \Hom_{K(\gth_+)} (\Hom (M, A^\bullet),
S^\bullet) = 0$, see \eqref{conj}, whenever $A^\bullet$ is
acyclic. This is true because the complex $\Hom(M, A^\bullet)$ is
acyclic and $S^\bullet$ is $K$-injective with respect to $\gth_+$. The
second is that $\Hom_{K(\gth_-)} (S^\bullet
\otimes M, B^\bullet) \linebreak[1] = \linebreak[0]
\Hom_{K(\gth_-)} ( S^\bullet, \linebreak[1] \Hom ( M, \linebreak[0]
B^\bullet)) = 0$ whenever $B^\bullet$ is an acyclic complex in
$K(\gth_-)$. That is true because the complex $\Hom(M, B^\bullet)= 0$
is acyclic as well, and $S^\bullet$ is $K$-projective as a complex of
$\gth_-$-modules. Thus the complex $(S^\bullet
\otimes M)^{\gth_+}_{\gth_-}$ computes the \si\ \coh\ $\sh (\gth, M)$.

Similarly, $S^\bullet \otimes \sindh M$ is a weakly \sj\ resolution of
$\sindh M$ over $\gtg$. Therefore the \coh\ of the complex $(S^\bullet
\otimes \sindh M)^{\gtg_+}_{\gtg_-}$ is equal to $\sh (\gtg, \linebreak[0]
\sindh M)$.

To conclude the proof, it suffices to observe the canonical
isomorphism
\[
(S^\bullet
\otimes M)^{\gth_+}_{\gth_-} = (S^\bullet
\otimes \sindh M)^{\gtg_+}_{\gtg_-},
\]
which comes from the functorial isomorphism
\[
(N \otimes M)^{\gth_+}_{\gth_-} = (N
\otimes \sindh M)^{\gtg_+}_{\gtg_-}
\]
for any $\gtg$-module $N$. This isomorphism is a universal property of
\si\ induction; it suffices to establish the following isomorphism of
right $\gtg$-modules:
\begin{equation}
\label{univ}
N \otimes^{\gtg_+}_{\gtg_-} \US = N
\end{equation}
as a $\gtg$-module. This statement is the computation of the
semi-invariants of the $\gtg$-module $N \otimes \US$, which follows
from Propositions~\ref{iso} and \ref{iso-1}. Indeed,
Proposition~\ref{iso} shows that $\US \cong U(\gtg_+)^* \otimes
U(\gtg_-)$ as a right $\gtg_-$-module. Therefore, taking the
$\gtg_-$-coinvariants of the tensor product with $N$ and $\lb$ will
reduce $N \otimes \US \otimes \lb$ to $N \otimes U(\gtg_+)^* \otimes
\lb$, which is canonically identified with $N \otimes U(\gtg_+)^*$ as
a vector space. Similarly, Proposition~\ref{iso-1} implies that the
$\gtg_+$-invariants of $N \otimes \US$ are computed as $N
\otimes U(\gtg_-)$. And it is obvious that the image of $N \otimes
U(\gtg_-)$ in $N \otimes U(\gtg_+)^*$ is exactly $N$. The fact that
\eqref{univ} is a $\gtg$-module isomorphism follows from the
naturality of the intermediate isomorphisms.
\end{proof}

\section{Wakimoto modules}

Wakimoto modules form an interesting class of highest-weight modules
over an affine Kac-Moody algebra, playing an intermediate role between
Verma and contragredient Verma modules.

\subsection{Contragredient Verma modules}

Before going into the subject of Wakimoto modules, we would like to
outline certain well-known properties of contragredient Verma modules,
mainly for motivational reasons. Let $\gtg$ be a simple
finite-dimensional Lie algebra and $\wgtg = \nc[z,z^{-1}] \otimes \gtg
\oplus \nc K \oplus \nc d$ the corresponding affine Kac-Moody algebra,
$K$ being its central element and $d= z\frac{\partial}{\partial z}$
(see V.~G. Kac \cite{kac:r} for more detail on affine Kac-Moody
algebras). Choose a Cartan subalgebra $\gtt$ in $\gtg$ and a system of
simple roots for $\gtg$. Let $\gtn_+$ ($\gtn_-$) be the subalgebra of
$\gtg$ spanned by the positive (respectively, negative) root
subspaces; then $\gtg = \gtn_+ \oplus \gtt \oplus
\gtn_-$. Furthermore, we can define a $\nz$-grading on $\wgtg$ by
putting $\deg 1 \otimes g = 1$ whenever $g$ is in a simple root space
of $\gtg$, $\deg z = \rank \gtg + 1$, $\deg 1 \otimes g = -1$ if $g$
is in the root space corresponding to a negative simple root, $\deg
z^{-1} = - \rank \gtg -1$, and $\wgtg_0 = \gtt \oplus \nc K \oplus \nc
d$. Then
\begin{align*}
\wgtg_+ & = \bigoplus_{n > 0} \wgtg_n = \gtn_+ \oplus (z\nc[z] 
\otimes \gtg) ,\\
\wgtg_- & = \bigoplus_{n \le 0} \wgtg_n = \gtn_- \oplus (z^{-1} \nc[z^{-1}]
\otimes \gtg) \oplus \gtt \oplus \nc K \oplus \nc d.
\end{align*}
The splitting $\wgtg = \wgtg_+ \oplus \wgtg_-$ along with the
functional $\beta$ such that $\beta (K) = 2h^\vee$, where $h^\vee$ is
the dual Coxeter number, $\beta|_{\gtt} = 2\rho$, $\beta(d) = 1$,
defines a \si\ structure on $\wgtg$.

Given a character $\lambda$ of $\wgtg_0$, one can define a
\emph{contragredient Verma module} $V^*(\lambda)$ using coinduction:
\[
V^*(\lambda) = \coin_{\wgtg_-}^{\wgtg} \nc_\lambda,
\]
where $\nc_\lambda$ is the corresponding one-dimensional
representation of $\wgtg_0$ extended by zero to a one-dimensional
representation of $\wgtg_-$. 
\begin{thm}
The module $V^*(\lambda)$ can be uniquely characterized as a
$\wgtg_0$-di\-ag\-o\-nal\-iz\-able $\wgtg$-module in the category
$\OO$, such that
\[
H^{i} (\wgtg_+, V^*(\lambda)) = \begin{cases}
	\nc_\lambda & \text{if $i =0$},\\
	0 & \text{otherwise}.
\end{cases}
\]
\end{thm}
\begin{proof}
The cohomological condition is satisfied for a contragredient Verma
module: it is obviously injective as a $\wgtg_+$-module.

Conversely, suppose we have a module $M$ satisfying the cohomological
condition. Let us prove it is isomorphic to the corresponding
contragredient Verma module.

Diagonalize the $\wgtg_0$-action and define a morphism $M \to
\nc_\lambda$ of $\wgtg_0$-modules which takes a highest-weight vector
$m_\lambda$ (one invariant under $\wgtg_+$), defined up to a scalar
factor, of $M$ to the generator of $\nc_\lambda$, mapping all the
other weight subspaces to zero. It will automatically be a morphism of
$\wgtg_-$-modules.

The constructed morphism induces a unique morphism $f: M \to
V^*(\lambda)$ of $\wgtg$-modules by virtue of the universality of
coinduction \eqref{universal}. The morphism $f$ must be injective,
otherwise it has a kernel, and we can consider the short exact
sequence
\[
0 \to \Ker f \to M \to \im f \to 0.
\]
Now let us look at the long exact sequence of \coh\ of $\wgtg_+$.
Since the highest-weight vector $m_\lambda$ of $M$ obviously maps
nontrivially to $\im f$, it follows that $H^0(\wgtg_+, \Ker f) = 0$
--- no highest-weight vectors, which implies $\Ker f = 0$, because
$\Ker f$ lies in the category $\OO$.

The morphism $f$ should be surjective, otherwise it has a cokernel:
\[
0 \to M \to V^*(\lambda) \to \Coker f \to 0.
\]
Again, the long exact sequence of \coh\ gives
\[
0 \to H^0(\wgtg_+, M) \to H^0(\wgtg_+, V^*(\lambda)) \to H^0(\wgtg_+,
\Coker f ) \to 0.
\]
Since $H^0(\wgtg_+, M) \to H^0(\wgtg_+, V^*(\lambda))$ is an
isomorphism, $\Coker f$ may not have a highest-weight vector, which
means it should be trivial.
\end{proof}

\begin{rem}
From the proof, we can observe that it is enough to require the
vanishing of only $H^1 (\wgtg_+, M)$ to make sure that $M \cong
V^*(\lambda)$.
\end{rem}

Analogously, one can homologically characterize the
\emph{Verma module} $V(\lambda) = \inn_{\wgtg_+ \oplus \wgtg_0}^\gtg
\nc_\lambda$ as a unique $\wgtg_0$-diagonalizable $\wgtg$-module in the
category $\OO$, such that
\[
H_{i} (\wgtg_{< 0}, V(\lambda)) = \begin{cases}
	\nc_\lambda & \text{if $i =0$},\\
	0 & \text{otherwise}.
\end{cases}
\]

\subsection{Wakimoto modules}

First, let us introduce a setup in which Wakimoto modules arise, see
Feigin-Frenkel \cite{ffr} for more detail. Consider an alternative
splitting $\wgtg = \gta \oplus \bar \gta$, where
\begin{align*}
\gta & = (\nc[z,z^{-1}] \otimes \gtn_- ) \oplus (z\nc[z] \otimes \gtt),\\
\bar \gta & = (\nc[z,z^{-1}] \otimes \gtn_+) \oplus 
(\nc[z^{-1}] \otimes \gtt) \oplus \nc K \oplus \nc d.
\end{align*}
One can think of $\bar \gta$ as a Borel subalgebra of $\wgtg$,
obtained as the limit under the action by the elements $2 m \rho$,
$\rho$ being the half-sum of positive roots of $\gtg$, $m
\to \infty$, of the affine Weyl group on the Borel subalgebra
$\wgtg_-$. The subalgebras $\gta$ and $\bar \gta$ are obviously
$\nz$-graded subalgebras of $\wgtg$, and the decomposition $\bar \gta
= \bar \gta_+
\oplus
\bar \gta_-$, where
\begin{align*}
\bar \gta_+ & = \bar \gta \cap \wgtg_+ = \nc[z] 
\otimes \gtn_+ ,\\
\bar \gta_- & = \bar \gta \cap \wgtg_- = (\nc[z^{-1}]
\otimes \gtn_+) \oplus (\nc[z^{-1}] \otimes \gtt) \oplus \nc K \oplus \nc d,
\end{align*}
along with the same functional $\beta$ induces a \si\ structure on the
subalgebra $\bar \gta$. Similarly, the decomposition $\gta = \gta_+
\oplus \gta_-$, where
\begin{align*}
\gta_+ & = \gta \cap \wgtg_+ = (z\nc[z] 
\otimes \gtn_-)  \oplus (z\nc[z] \otimes \gtt),\\
\gta_- & = \gta \cap \wgtg_- = \nc[z^{-1}]
\otimes \gtn_-,
\end{align*}
with $\beta = 0$ defines a \si\ structure on the subalgebra $\gta$.

Now we are ready to give a constructive definition of a Wakimoto
module.
\begin{df}
Let $\lambda$ be a character of $\wgtg_0$ and $\nc_\lambda$ the
corresponding one-di\-men\-sion\-al $\wgtg_0$-module, extended
trivially to $\bar \gta$. Then the \emph{Wakimoto module} $W(\lambda)$
is the semi-induced module
\[
W(\lambda) = \sind_{\bar \gta}^{\gtg} \nc_\lambda.
\]
\end{df}

\subsection{The structure of a Wakimoto module}

Let $\alpha_1, \alpha_2, \dots$ be all the roots of the Lie algebra
$\wgtg$ and $e_{\alpha_s,i_s}$ a basis of the root subspace
$\wgtg_{\alpha_s}$, $i_s = 1, \dots, m_s$, where $m_s$ is the
multiplicity $\mult \alpha_s$, which is one when the root $\alpha_s$
is real, $\rank \gtg$ when $\alpha_s$ is imaginary, and $\rank \gtg +
2$ when $\alpha_s = 0$. The Poincar\`e-Birkhoff-Witt Theorem states
that the products
\[
e^{n_{1,1}}_{\alpha_1,1} \dots e^{n_{1,m_1}}_{\alpha_1,m_1}
e^{n_{2,1}}_{\alpha_2,1} \dots e^{n_{2,m_2}}_{\alpha_2,m_2} \dots
\]
containing only a finite number of terms form a basis of the \uea\
$U(\wgtg)$. Let
\[
(e^{n_{1,1}}_{\alpha_1,1} \dots e^{n_{1,m_1}}_{\alpha_1,m_1}
e^{n_{2,1}}_{\alpha_2,1} \dots e^{n_{2,m_2}}_{\alpha_2,m_2} \dots)^*
\]
denote the corresponding element of the dual basis in
$U(\wgtg)^*$. Let $\beta_1, \beta_2, \dots$ be the roots of $\gta_+$
with the corresponding multiplicities $k_1, k_2, \dots$ and $\gamma_1,
\gamma_2, \dots$ the roots of $\gta_-$. The roots of $\gta_-$ are all
real, therefore, the multiplicities are all one. Let $w_\lambda$ be
the generator of $\nc_\lambda$.

\begin{prop}
\label{structure}
\begin{enumerate}
  \item
  \label{1}
In the notation of the previous paragraph, the vectors
\begin{equation}
\label{basis}
(e^{n_{1,1}}_{\beta_1,1} \dots e^{n_{1,k_1}}_{\beta_1,k_1}
e^{n_{2,1}}_{\beta_2,1} \dots e^{n_{2,k_2}}_{\beta_2,k_2} \dots)^*
\otimes (e^{l_{1}}_{\gamma_1} \dots e^{l_2}_{\gamma_2} \dots )
\otimes w_\lambda
\end{equation}
form a basis of the Wakimoto module $W(\lambda)$.
  \item 
The \emph{formal character} $\ch W(\lambda) = \sum_{\mu \in \wgtg_0^*}
\dim W(\lambda)_\mu e^\mu$, where $W(\lambda)_\mu$ is the weight space, of
$W(\lambda)$ is equal to
\[
\ch W(\lambda) = e^\lambda \prod_{\alpha \in \Delta_+}
(1-e^{-\alpha})^{-\mult \alpha},
\]
where $\Delta_+$ are the positive roots of $\wgtg$. Thus, $\ch
W(\lambda) = \ch V(\lambda) = \ch V^*(\lambda)$.
\end{enumerate}
\end{prop}

\begin{proof}
Since $U(\wgtg) = U(\gta_+) \otimes U(\gta_-) \otimes U(\bar \gta_+)
\otimes U(\bar \gta_-)$ as a vector space, a Wakimoto module
$W(\lambda)$ is isomorphic to $U(\gta_+)^* \otimes U(\gta_-) \otimes
\nc_\lambda$, and the Poincar\'e-Birkhoff-Witt Theorem for $\gta_+$ and
$\gta_-$ proves \ref{1}.

The weight of a vector \eqref{basis} is obviously $- ((n_{1,1} + \dots
+ n_{1,k_1}) \beta_1 + (n_{2,1} + \dots + n_{2,k_2}) \beta_2 + \dots)
+ (l_1 \gamma_1 + l_2 \gamma_2 + \dots) + \lambda$. Notice that
$-\beta_1, -\beta_2, \dots, \gamma_1, \linebreak[0] \gamma_2, \dots$
comprise all the negative roots of $\wgtg$. Therefore,
\[
\ch W(\lambda) = e^\lambda \prod_{\alpha \in \Delta_+}
(1+e^{-\alpha}+e^{-2\alpha} + \dots)^{\mult \alpha},
\]
which implies the character formula.
\end{proof}

\subsection{Semi-infinite \coh\ of Wakimoto modules}

\begin{thm}
The Wakimoto module $W(\lambda)$ is a
$\wgtg_0$-di\-ag\-o\-nal\-iz\-able $\wgtg$-module from the category
$\OO$, such that
\[
H^{\infty/2 + i} (\gta, W(\lambda)) = \begin{cases}
	\nc_\lambda & \text{if $i =0$},\\
	0 & \text{otherwise}.
\end{cases}
\]
\end{thm}

\begin{proof}
It follows from Proposition~\ref{structure}, that as an $\gta$-module,
$W(\lambda) \linebreak [2] = \linebreak[1] \US(\wgtg)
\linebreak[0]\otimes^{\bar \gta_+}_{\bar \gta_-} \nc_\lambda \cong
\US(\gta) \otimes \nc_\lambda$. By the \si\ Shapiro
Lemma applied to $\gtg = \gta$ and $\gth = 0$, we have the required
computation of the \si\ \coh\ of $W(\lambda)$.
\end{proof}

This Theorem was used by Feigin and Frenkel \cite{ffr} as the
definition of a Wakimoto module, without proving that that
cohomological property defined it uniquely. This uniqueness of a
Wakimoto module is presumably true, but it would be more cautious to
call it an open problem. In a recent e-mail message to the author,
S.~M. Arkhipov has suggested an outline of a prospective solution,
using the fact that a Wakimoto module is a projective limit of twisted
contragredient Verma modules.


\providecommand{\bysame}{\leavevmode\hbox to3em{\hrulefill}\thinspace}

\end{document}